\documentstyle[12pt]{article}
\begin{document}

\title{Quasi-continuous symmetries of non-Lie type}

\author{Andrei Ludu and Walter Greiner  \\
\normalsize{Institut f\"ur Theoretische Physik der } \\ 
\normalsize{J. W. Goethe - Universit\"at, 
Robert-Mayer-Strase 8-10, D-60054, Frankfurt am Main,} \\
\normalsize{ Germany}
}

\date{}
\maketitle
\begin{abstract}
We introduce a smooth mapping of some
discrete space-time symmetries into quasi-continuous ones.
Such transformations are related with q-deformations
of the dilations of the Euclidean space and with the non-commutative space. 
We work out two  examples of Hamiltonian invariance under such symmetries. 
The Schr\"odinger equation for a free particle is investigated in such a 
non-commutative plane and a connection with anyonic statistics is found.

\end{abstract}

\setlength{\baselineskip} {4ex}.
 
\vskip2cm

\vfill
\eject

Continuous symmetries are generally described in terms of Lie
groups and algebras through their irreducible unitary representations acting
on the Hilbert space of states. Among these, the space-time symmetries are 
expressed as the 
invariance of the dynamical equations under transformations of
the coordinate frame.
One advantage of continuous
symmetries over the discrete ones appears clearly
in  the Lagrangean field theories where, due to the Noether
theorem, one can define currents and conserved quantities 
associated with  the time evolution
of the dynamical system.
The discrete symmetries arising in physiscs can be clasificated in $Z_{2}$
graded symmetries consisting in inversions of the space (parity P,
reflexions, mirroring), charge conjugation (C) and time reversal (T) on one
side, and permutations and braids on the other side. The former act on 
one-particle states or directly on the space-time manifold. The operators
associated with these symmetries have the square equal with one.
The latter ones act on many-particle states, and they are connected with the 
statistics of the physical system; they are finite or infinite dimensional
unitary irreducible representations of some discrete groups over some Hilbert 
space of states.
  
In this paper we introduce a continuous algebra of transformations 
in which are embedded some discrete space-time symmetries. 
A first direct way to construct such a structure is
to introduce an associative 
continuous algebra of operators, defined by some 
 commutator relations, containing the identity and
the discrete transformations among its elements. 
Such an algebraic-continuous structure is more general then Lie
algebras or loop groups. 
If both the identity and the discrete symmetries could be brought together
in the same smooth algebraic structure, it is possible that new symmetries 
(associated with the intermediate steps in between the 
identity-discrete limits)  could occure.
We show in the following that such a structure exists and can be obtained
as a q-deformed algebra. 
Quantum groups [1-3] (quantized universal enveloping 
algebras, q-algebras, q-deformations)
have been the subject of numerous recent studies in mathematics and physics.
They represent some
special deformations ($q\neq 1$) of the universal
enveloping algebra of Lie algebras ($q=1$) [1-3]. 
For a recent review see [4]. 
In the present paper
we are both interested  in some limiting cases
of the q-deformations, which should 
meet the discrete original symmetries, and in the intermediate
(nongeometric) symmetries (joining the discrete elements).
The mapping of the discrete symmetries on the  continuous
ones  
results in interesting consequences 
for the space-time discretisation, non-commutative spaces and
spontaneous breaking of symmetry.

We introduce 
the $R_{3}$ Euclidean (commutative) space
generated by $(x_{i})=(x,y,z)$ and its tangent space, generated by
$\partial _i=\partial /\partial x_i$. 
Twelve infinitesimal (Lie) generators 
${\hat {\xi}} \in \{ {\partial }_{i}$, $x_{j}{\partial }_{i} \} $, 
 act on the product space of $R_3$ with its tangent space.
With these generators one 
constructs different Lie algebras with action on the functions defined on
$R_{3}$.
We associate to each generator ${\hat {\xi }}$ 
a one-real parameter Lie group, 
$g(\phi )|_{\phi \in R} $ by way of the exponential map ${\hat {\xi }}
\rightarrow g(\phi )=exp(\phi {\hat {\xi }})$.
In order to play an example, we restrict to the following 
discrete transformations, acting on $R_{2}$:
identity $1$, infinitesimal generator of rotations in plane $v$
(or finite rotation
of $-\pi /2$), reflection against $Oy$ axis $r_y$ and permutation $P$
(or reflection
against the line $x=y$). They satisfy the relations: $r_{y}^{2}=P^{2}=1$,
$v^{2}=-1$, and (excepting $I$) they  
generate a Lie algebra, isomorphic with
$su(2)$, and defined by the commutators:
\begin{eqnarray}
[v,r_{y}]=2P, \,\,\, [P,v]=2r_{y}, \,\,\, [P, r_{y}]=2v. \hfill
\end{eqnarray}
The isolated (discrete symmetries) $P$ and $r_{y}$, can not
be regarded as group elements, since they do not belong to any representation
of $SU(2)$, hence they are not realised in terms of a certain
$g(\phi )$, for any $\phi $. 
Consequently, following the Lie
approach,
one can not smootly map the identity group element into $P$ or $r_{y}$.
In fact, this conclusion is the algebraic formulation of the geometric
imposibility of
continuously mapping the $y$ axis into the $-y$ axis
through a Lie transformations (rotations against the $x$ axis, dilation of
the $y$ axis, etc), in two dimensions only.
All such methods always yield to intermediate situations having
algebraic singularities: the number of dimensions increases to three or
decreases to one.
From the topological point of view
the mirror transformation $r_{y}$ 
is similar with the problem of mirroring a knot [5].
The right-handed and the left-handed coordinate frames of the plane are 
actualy oriented knots since one of the axis always
undercrosses the other. The coordinate frame is an oriented
 knot which holds the same topological informations as the original image.
So, the discrete transformation $r_{y}$ of the coordinate frame, can be
realised as a mirror transformation of a knot (the interchanging of the
roles of undercrossed/overcrossed at each knot of a diagram).
There are
examples of (achiral) knots which can be continuously deformed into their 
mirror image (Eight knot) and examples of (chiral)
knots which can not (Trifoil knot) [5].
This continuous deformation is a finite succesion of
Reidemeister moves of ambiental
isotopy (transformations in the plane which simulate the corresponding
natural topological transformations in the 3-dimensional space of the
unfolded correspondent of the knot).
Consequently, an algebraic treatment of the Reidemeister moves
can direct the problem to the construction of a
continuous connection of the discrete mirrorings with the identity.
Since the quantum group $gl_q (2)$  
provides representations of such moves [5], we construct in the following
special q-deformations of the Lie algebra $gl(2)$ in relation 
with knots and braids.

In order to use such deformations we introduce the q-deformed commutators in
the form $[x,y]_{q}\equiv xy-qyx$, $[x,y]_{q=1}=[x,y]$, [4]. Such relations
exist in the universal enveloping algebra of all the products
of the generators, and they generate a special
type of quantum group.
By introducing the transformations.
$R_{y}=r_{y}+v$, $V=r_{y}-v$ we have the new commutators:
\begin{eqnarray}
[R_{y},P]_{q}=(1+q)R_{y}, \,\,\, [V,P]_{q}=-(1+q)V, \,\,\, 
[R_{y},V]_{q}=2((1-q)1-(1+q)P). \hfill 
\end{eqnarray}
We remark that in the last commutator we have obtained a deformed element
which maps  $1$ into $P$, when $q$ goes from $q_{1}=-1$
to $q_{2}=1$. 
An example of action of the above defined deformed commutators
on the real line is given by introducing
a quantum analogue of the $x$ coordinate through the operator
\begin{eqnarray}
{\hat x}=x{\hat {Q}}(x\partial _x,q),
\end{eqnarray}
where ${\hat {Q}}(x \partial _x, q)$ is a linear operator given 
by a function of the infinitesimal 
generator of dilations $x\partial _x$ and deppending on the parameter $q$,
such that $\lim_{q \rightarrow 1}{\hat {Q}}=1$. Following [6,7]
the corresponding canonical momentum is defined as
\begin{eqnarray}
{\hat {\partial }}_x={\hat {Q}}(x \partial _x,q)\partial _x.
\end{eqnarray}
From eqs.(3,4) it results the following commutation relation:
\begin{eqnarray}
[{\hat {\partial }}_x,{\hat x}]=1+(q^2-1){\hat x}{\hat {\partial }}_x,
\end{eqnarray}
or, written in the formalism of q-deformed commutators,
\begin{eqnarray}
[{\hat {\partial }}_x,{\hat x}]_q =1.
\end{eqnarray}
Eq.(6) becomes the usual commutator relations between
the $x$ and $\partial _x$ operators, in the limit $q\rightarrow 1$.
From eq.(6) the operator $\hat {Q}$ defined in eq.(4) must satisfy the
condition
\begin{eqnarray}
{\hat {Q}}^2 +{\hat {Q}}x\partial _x {\hat {Q}}=1+q^2 x {\hat {Q}}^2 
\partial _x,
\end{eqnarray}
which, applied on integer real functions $f(x)=\sum_{j}f_j x^j$, reads
\begin{eqnarray}
(j+1){\hat {Q}}^2 (j)=1+q^2 j{\hat {Q}}^2 (j-1),
\end{eqnarray}
From eq.(8) we obtain the solution
\begin{eqnarray}
{\hat {Q}}^2 (j)=q^j {{[j+1]} \over {j+1}},
\end{eqnarray}
where $[j+1]$ is the q-deformation of $j+1$, [3,4].
Consequently, we have the realisation of 
the new coordinate operator in the form
\begin{eqnarray}
{\hat x}=x
\sqrt
{q^{x\partial _x}{{q^{2x\partial _x}-1} \over {(q-1)(x\partial _x +1)}}}.
\end{eqnarray}
similarly with the realisation introduced in [7].
In order to evoid the square root 
occuring in eq.(10)
we can modify the q-deformed commutator relations, preserving the same 
behaviour in the limit $q=1$.
If the new coordinates of the phase-space satisfy
\begin{eqnarray}
[{\hat {\partial }}_x, {\hat { x}}]=q^{2{\hat x}{\hat {\partial }}_x},
\end{eqnarray}
by following the same procedure similar with eqs.(6,7) we get for 
${\hat {Q}}$ a simpler form:
\begin{eqnarray}
{\hat {Q}} =q^{x\partial _x},
\end{eqnarray}
which gives us just the dilation operator.

In order to generalise the above construction 
we have to introduce a 
set of continuous transformations of the coordinates, depending 
on one complex parameter $q$, having no Lie algebraic equivalent 
and which, for certain fixed values of $q$,
approache the space inversions 
${\hat {I}}_x$, ${\hat {I}}_y$, ${\hat {I}}_z$.
These general structures are nonlocal, non-Lie and non-linear transformations,
 i.e. they do not form a Lie group (they are not of the form
$e^{\epsilon v}$ with $v$ a vector field of a Lie algebra not depending on
$\epsilon $).
 
We introduce over $R_3$ the
infinitesimal generators of the dilation ${\hat {d}}_i=x_i \partial _{x_i}$
and their corresponding one-dimensional
Lie groups ${\hat {D}}_k(s)=q^{{\hat {d}}_k}=e^{is {\hat {d}}_k}$, $[{\hat {D}}_k,{\hat {D}}_j]=0$,
where we use a complex
deformation in the form $q=e^{is}$ with $s\in [0,\pi ]$.
The action of the dilation operators on analytical functions on $R_3$
is continuous with respect to $s$.
In the limit $q\rightarrow 1$ ($s \rightarrow 0$) we have ${\hat {D}}_{i}
\rightarrow 1$ and in the limit $q\rightarrow -1$ ($s \rightarrow \pi$)
we have ${\hat {D}}_i \rightarrow {\hat {I}}_{i}$.
Consequently, in the range $s:[0 , \pi ]$
the operators ${\hat {D}}_{i}$ smoothly map  the unit element into the 
corresponding inversion operator ${\hat {I}}_i$, 
${\hat {D}}_{i}(\pi )={\hat {I}}_i$ generating in this way all the
inversions of $R^{3}$.
A q-deformed operator integer function ${\hat {Q}}({\hat {d}}_i,s)=\sum _{k}C_k (s){\hat {d}}_{i}^{k}$ 
has the action
on real integer functions $f(x_i)=\sum _{j}f_j x_i^j$,
${\hat {Q}}({\hat {d}}_i,s)f(x_i)=\sum _{j}
f_j{\hat {Q}}(j,s)x_i^j $ which
approaches $f(x_i)$ for ${\hat {Q}}({\hat {d}}_i,1)
= 1$.
This operator should be linear and should have the following limits:
${\hat {Q}}({\hat {d}}_i,1)=1, {\hat {Q}}({\hat {d}}_i,-1)=
e^{i\pi {\hat {d}}_i}={\hat {I}}_i$.

In the following we we want to find the most general one-dimensional 
Hamiltonian having a symmetry for these types of transformations. 
In the one-dimensional case ($x_i=x$)
a  Hamiltonian ${\hat {H}}$ 
is invariant to the symmetry introduced by the operator
${\hat {Q}}={\hat {Q}}({\hat {d}}_i,s)$ if $[{\hat {H}},{\hat {Q}}]=0$. 
We choose a general one-dimensional Hamiltonian in the form
\begin{eqnarray}
{\hat {H}}=-\partial _{x}^2+V(x)+W(x \partial _{x}),
\end{eqnarray}
and the corresponding Schr\"odinger equation ${\hat {H}}f(x)=Ef(x)$.
The Hamiltonian consists of the sum of the kinetic energy term, a
local potential $V(x)=\sum_{j}V_j x^j$
and an effective potential $W(x \partial _{x})=\sum W_j (x \partial _{x})^j$.
The last term  commutes  with ${\hat {Q}}$. 

The condition of symmetry of $\hat H$ under the transformation ${\hat {Q}}$,
$[{\hat {Q}},{\hat {H}}]f(x)=0$ applied to an arbitrary integer function
$f(x)=\sum f_k x^k$, after identification of
the corresponding powers of $x$, reads
\begin{eqnarray}
f_{k+2}(k+1)(k+2)({\hat {Q}}(k)-{\hat {Q}}(k+2))=\sum _{j=0}^{k-1}V_{k-j}f_j({\hat {Q}}(k)-{\hat {Q}}(j)).
\end{eqnarray}
Eq.(14) must be solved together with the Schr\"odinger equation for $f$
which, expanded in powers of $x$, reads:
\begin{eqnarray}
f_{k+2}(k+1)(k+2)=\sum_{l=0}^{k}f_l V_{k-l}-Ef_k +f_k W(k).
\end{eqnarray}
In the simplest case ${\hat {Q}}={\hat {D}}_x$ eq.(14) becomes
\begin{eqnarray}
f_{k+2}={{\sum_{j=0}^{k-1}V_{k-j}f_j (1-q^{j-k})} \over {(1-q^2)(k+1)(k+2)}}.
\end{eqnarray}
Eqs.(14,15) give the conditions for the existence of the continuous symmetry
described by the operator ${\hat {Q}}$. Both these equations are identities for ${\hat {Q}}=1$.
In the case of the inversion $x \rightarrow -x$, ${\hat {Q}}\rightarrow {\hat {I}}_x$, 
i.e. ${\hat {Q}}(k)=(-1)^k$,
eq.(14) asks for $V_k=f_k=0$, for $k$ odd, ($V(-x)=V(x),f(x)=f(-x)$) like in
the traditional case.
In the free particle case we have $V(x)=0$ and from eq.(14) 
we get
\begin{eqnarray}
{\hat {Q}}(k)={\hat {Q}}(k+2),
\end{eqnarray}
and ${\hat {Q}}(k_{odd})=\pm 1$
and ${\hat {Q}}(k_{even})=1$. This restricts the allowed symmetries for the free
particle to inversion only. 
For
any potential of the form $V(x)=x^n$ we obtain the corresponding
invariance
condition $q^n=1$.
In this case the number of admissible discrete 
symmetries between $1$ and ${\hat {I}}_x$ is finite ($n$ roots of the unity)
and still there is no way of continuous mapping of the identity into the
mirroring. 

We want to solve eqs.(14-15) in the simple case of 
${\hat {Q}}={\hat {D}}$. From eq.(14) we note that a full invariance of the
wavefunction $f(x)$ occures if the potential $V(x)$ depends also on $q$.
Hence, we introduce and arbitrary
potential $V^0(x)=\sum_k V^0_k x^k$ independent of $q$,
 and we write the  potential
of the Hamiltonian in eq.(13) in the form of a transformation 
\begin{eqnarray}
V(x) \rightarrow V(x,q)=\sum_k V_k(q)x^k=
{1 \over 2}\sum_{k}kV^0_k {{q^2 -1} \over {1-q^{-k}}}x^k, 
\end{eqnarray}
where the coefficients $V^0_k$ (independent of $q$) determine the limit 
$V(x,1)=\sum_k V^0_k x^k=V^0 (x)$. Due to this choice, eq.(14) 
does not depend any more
on $q$ and one can obtain the coefficients $f_k$ as function of $V^0_k$ only.
This procedure makes sense only if the coefficients $1-q^{k-j}$ in eq.(16)
are not zero. This condition depends on $j,k$ and consequently on the nonzero
coefficients of $V(x,q)$ in its Taylor expansion, and on the values of $q$.
If we do not have such singularities (e.g. for $s$ 
irrational multiples of $\pi $)
the corresponding wavefunctions do not depend on $q$, i.e. do not depend on 
the global transformation
${\hat {Q}}$, for any $q$, from the identity ($q=1$) to the inversion ($q=-1$):
\begin{eqnarray}
f_{k+2}={{\sum_k V^0_{k-j}f_j} \over {(k+1)(k+2)}}.
\end{eqnarray}
The physical system remains
in the same quantum state under the action of all ${\hat {Q}}'s$ from 
${\hat {1}}$ to ${\hat {I}}_x$.
With the coefficients $f_k$ 
of the wavefunction obtained from eq.(19) we can determine 
the coefficients $W(k)$ and $E$ in eq.(15).
Thus, we obtain a dynamical symmetry for the Hamiltonian in eq.(13)
under ${\hat {Q}}(x\partial _x,q)$,
with its eigenfunctions $f(x)$ and eigenvalues $E$ independent of $q$.
As a consequence, the 
transformations ${\hat {Q}}(x\partial _x,q)$ keep invariant the physical states.
However the associated potential transforms with $q$, like in the case of a 
gauge transformation.
The transformations of the potential, given by eq.(18), can be 
written explicitely in the form of a nonlocal operator applied on $V^0 (x)$,
for $|q|<1$:
\begin{eqnarray}
V(x,q^2 )={{q^2 (q+q^{-1})} \over {2}}q^{x\partial _x }
\int {{dV^0 (x)} \over {dx}}d_q x,
\end{eqnarray}
that is the q-primitive of the derivative of $V^0$,
evaluated in $q^{1/2}x$, [8]:
\begin{eqnarray}
\int f(x)d_q x = (q^{-1}-q)x\sum_{n=0}^{\infty }q^{2n+1}f(q^{2n+1}x).
\end{eqnarray}
For $q=1$ the q-primitive tends to the normal integration
and eq.(20) describes the action of the identity operator on $V$.
This transformation of the 
potential is a sort of a nilpotent operation:
one acts first with an operator (the derivative and a scaling in $x$ with
the factor $q^{1/2}$) and then with a q-deformation
of the inverse of this operator (q-integration). The result is not the indentity
but a sort of a "defect" of the identity: an infinitesimal derivative followed by
a finite-difference integration and scaling. 
All these results obtained for $V(x)$ and $W$ 
can be generalised to three dimensions, too.
For $q\in R$ the effect of the transformation introduced in eq.(29)
is a scalling in $x$. For potentials having a pole in $x_0 $ the
transformation moves the pole in $x_0 /q^{1/2}$ and for $q\rightarrow 0$
the pole is eliminated. For $q\in C$ the transformed potential becomes
complex and the poles are translated into the imaginary
extension of the $x$ axis.
We present such an example, for a Coulomb-like
potential, for real deformations in Fig.1 and  for complex
deformations in Figs.2. In the complex case 
($q=e^{is}$, $s\in R$) the real part of the
q-deformed (${\hat {Q}}$-invariant) potential behaves completely different
from the original Coulomb potential, for $q$ not a root of the unity,
and transforms into a bounded potential. In the limit $q=-1$, $s=\pi $ 
the pole is translated from $1$ to $-1$, as the general formalism asks. 

We analyse now the  situation when $s$ goes from 0 to $\pi $
and provides zeros for $q^{2(k-j)}-1$ in eq.(16). 
We introduce the concept of quasi-continuous
transition between the discrete symmetries.
Let us define a discrete equidistant partition of 
the interval $[0,\pi ]$, $\triangle _N$,
by the points $s_{n,N}={n \over N}\pi $, $n=0,1,...,N$ and  consider that
$s$ goes from 0 to $\pi $ taking only the values $s_{n,N}$.
In this case, for any partition $\triangle _N$, we can define a complete
invariant potential $V(x,s_{n,N})$ denoted $V_N(x)$
\begin{eqnarray}
V_N (x)=\sum_{j=0}^{\infty }A_j x^{2jN}+
\sum_{j=0}^{\infty }B_j x^{(4j+1)N}+\sum_{j=0}^{\infty }C_j x^{(4j+3)N},
\end{eqnarray}
with arbitrary coefficients $A_j ,B_j$ and $C_j$.
The coefficients $A_j$ do not contribute in eq.(14), since 
$[(j-k)/2]_{s=s_{n,N}}=0$. In the same way we note that the coefficients
$B_j$, $C_j$ occure in eqs.(14,16) as $1$ and $-1$, respectively.
Hence, we can write,
in the partition $\triangle _N$, eqs.(14,16) in the form:
\begin{eqnarray}
f_{k+2}={
{\sum_{j=0}^{\biggl [ {{k-3} \over {4}}\biggr ] }
\biggl (B_{j}f_{k-4j-1} -C_j f_{k-4j-3} \biggr )
} \over 
{e^{is_{n,N}}sin(s_{n,N})(k+1)(k+2)}},
\end{eqnarray} 
where the right brakets in the limit of the sum represent the integer part.
In each partition $\triangle _N$, eqs.(14,16) allow the obtaining of the eigenfunctions
from the potential coefficients only, without any dependence on $s$.
By introducing the obtained  coefficients $f_j$ in eq.(15) 
we can solve all the system,
i.e. we deduce the $W(j)$'s. 
In this case we have an exact quasi-continuous symmetry for the transformation
${\hat {Q}}$ in the partition $\triangle _N$. It is exact
because all functions, the Hamiltonian and the
wavefunctions, do not depend on any value of $s$ in the given partition,
and it is quasi-continuous because $s$ takes only discrete values.
For higher values of $N$ this symmetry becomes very close to a continuous one.

By generalising eqs.(3,4,10) to
the continuous operators ${\hat {x}}_i=x_i {\hat {Q}}_i({\hat {d}}_i,q)$
and ${\hat {p}}_i={\hat {Q}}_i({\hat {d}}_i,q)\partial _i$
we have:
$$
{\hat {Q}}_{x}={\hat {Q}}({\hat {d}}_x,s)q^{{\hat {d}}_{y}+{\hat {d}}_{z}}, 
\ \ \
{\hat {Q}}_{y}={\hat {Q}}({\hat {d}}_y,s)q^{{\hat {d}}_{z}}, \ \ \ 
{\hat {Q}}_{z}={\hat {Q}}({\hat {d}}_z,s), \hfill
$$
\begin{eqnarray} 
{\hat {Q}}({\hat {d}}_i,s)=
exp({{s{\hat {d}}_i } \over 2})[{\hat {d}}_i+1]^{1/2} 
({\hat {d}}_i +1)^{-1/2},
\end{eqnarray}
where $[x]={{\sin (sx)} \over {\sin (s)}}$ represents the q-deformation
of the object $x$ (a c-number, an operator, etc.), [3-8].
In the limit $q\rightarrow 1$ ($s \rightarrow 0$) we have ${\hat {Q}}_{i}
\rightarrow 1$ and in the limit $q\rightarrow -1$ ($s \rightarrow \pi$)
we have
\begin{eqnarray}
{\hat {Q}}_{x} \rightarrow {\hat {I}}_{y}{\hat {I}}_{z}=-{\hat {I}}_{x}, \,\, 
{\hat {Q}}_{y} \rightarrow {\hat {I}}_{z}, \,\, 
{\hat {Q}}_{z} \rightarrow 1.\hfill
\end{eqnarray}
The generators ${\hat {Q}}_{i}$ form an associative commutative
algebra, $[{\hat {Q}}_{i},{\hat {Q}}_{j}]=0$. 
Consequently, in the range $s:[0 , \pi] $
the operators ${\hat {Q}}_{i}$ map smoothly the unit element into a certain
reflection operator and ${\hat {Q}}_{i}(\pi )$ generate in this way, together with
the rotations of the space, all possible
inversions in $R^{3}$.

A first 
interpretation for the operators given in eq.(24) 
comes out from traditional quantum mechanics. To the classical coordinate
observables $(x,y,z)$ one associates some operators, not like in the
traditional way (multiplicative coordinate functions ${\hat x}_{i}=
x_{i}{\hat 1}$) but through the q-deformed operators ${\hat {x_{i}}}=
{x_i}{\hat {Q}}_{i}$. 
The action of ${\hat {Q}}({\hat {d}}_i ,s)$ on integer functions of $x_i $
is given by ${\hat {Q}}({\hat {d}}_i ,s)x_{i}^{k}=Q(k,s)x_{i}^{k}$.
Using also the fact that 
$exp(is{\hat {d}}_{i})x_j =exp(is)x_{j}\delta _{i,j}$,
we obtain the action
$$
{\hat {Q}}_{i}x_{j}=\biggl (
\sqrt{{{q^2 +1} \over {2}}}\delta _{i,j}+q(\delta _{i,j-1}+\delta _{i,j-2})
+\delta _{i,j+1}+\delta _{i,j+2}
\biggr )x_{j},
$$
for ${\hat {Q}}_{i}={\hat {Q}}_{x}, \dots $.
For some limiting cases we get exactly the discrete mirror operators:
${\hat {Q}}_{x}(x,y,z)|_{s=\pi }\rightarrow (x,-y,-z)$,
${\hat {Q}}_{x}(x,y,z)|_{s=\pi /2 }\rightarrow (0,iy,iz)$,
${\hat {Q}}_{y}(x,y,z)|_{s=\pi }\rightarrow (x,y,-z)$,
${\hat {Q}}_{y}(x,y,z)|_{s=\pi /2}$ 
$\rightarrow (x,0,iz)$,
${\hat {Q}}_{z}(x,y,z)|_{s=\pi }\rightarrow (x,y,z)$,
${\hat {Q}}_{z}(x,y,z)|_{s=\pi /2}\rightarrow (x,y,0)$.
We can see that these operators behave also like projectors, cancelling some
space components, or like analytical prolongations of the real coordinates
into the complex plane. 
A second interpretation is that the operators ${\hat x}_{i}$ 
are the new coordinates of a q-deformed non-commutative space.
In this case  [6,7,9] we introduce the relations 
\begin{eqnarray}
{\hat x}_{i}{\hat x}_{j}=q{\hat x}_{j}{\hat x}_{i}, 
 \,\,\, i<j   \hfill
\end{eqnarray}
which generate a space in these new non-commutative coordinates.
We can express these new coordinates in the limit $s \rightarrow \pi$:
${\hat x} \rightarrow {\hat {I}}_{x}x$, ${\hat y} \rightarrow {\hat {I}}_{z}y$
and ${\hat z}\rightarrow z$. Similarly one may introduce the
momentum operators associate with the
non-commutative coordinates ${\hat x}_{i}$
in the form
\begin{eqnarray}
{\hat {\partial }}_{x}=q^{{\hat {d}}_{y}+{\hat {d}}_{z}}{\hat {Q}}({\hat {d}}_x,s){\partial }_{x},
\ \ \ \ 
{\hat {\partial }}_{y}=q^{{\hat {d}}_{z}}{\hat {Q}}({\hat {d}}_y,s){\partial }_{y},
\ \ \ \ 
{\hat {\partial }}_{z}&=&{\hat {Q}}({\hat {d}}_z,s){\partial }_{z}. \hfill
\end{eqnarray}
In the limiting cases we have
${\hat {\partial }}_{i}(s\rightarrow 0)={\partial }_{i}$ and
${\hat {\partial }}_{x}(s\rightarrow \pi )=-{\hat {I}}_{x}{\partial }_{x}$,
${\hat {\partial }}_{y}(s\rightarrow \pi )={\hat {I}}_{z}{\partial }_{y}$,
${\hat {\partial }}_{z}(s\rightarrow \pi )={\partial }_{z}$.
Consequently, the new coordinates and the corresponding derivatives obey the
non-commutative calculus [6,9]:
$${\hat {\partial }}_{i}{\hat x}_{j}=q{\hat {\partial }}_{j}
{\hat x}_{i}
$$
$${\hat {\partial }}_{i}{\hat x}_{i}-q^{2}
{\hat x}_{i}{\hat {\partial }}_{i}=1+(q^{2}-1)\sum_{j>i}{\hat x}_{j}
{\hat {\partial }}_{j}
$$
\begin{eqnarray}
{\hat {\partial }}_{i}
{\hat {\partial }}_{j}&=&q^{-1}
{\hat {\partial }}_{j}
{\hat {\partial }}_{i},
\end{eqnarray}
with $i\neq j$. 
Both eqs.(3) and (27) are invertible with respect to the map
${\hat {x}}_{i}(x_{j},\partial _{k})$,
${\hat {\partial }}_{i}(x_{j},\partial _{k})$ 
$\leftrightarrow $
$x_{i}({\hat {x}}_{j},{\hat {\partial }} _{k}),
\partial _{i}({\hat {x}}_{j}$, ${\hat {\partial }}_{k})
$, [7].
The action of the operators ${\hat {Q}}({\hat {d}}_i,s)$ on integer functions 
$f(x_i )=\sum_j f_j {x}^{j}_{i}$ is given by:
\begin{eqnarray}
{\hat {Q}}({\hat {d}}_i,s)f(x_i)=
\sum_j f_j{\hat {Q}}(j,s){x_i}^j_i={{F(q^2x_i)-F(x_i)} \over 
{q(q-q^{-1})x}}={\hat {D}}_q{\biggl ( \int f(x_i)dx_i\biggr ) }_{q^{-1}x},
\end{eqnarray}
where $F(x_i)=\int f(x_i)dx_i$ is the primitive of $f$,
and the operator ${\hat {D}}_qf(x)={{f(qx)-f(q^{-1}x)} \over {(q-q^{-1})x}}$
is the q-derivative [3-8] and reduces 
to the normal derivative in the limit
$q=1,s=0$.
In the limit $s\rightarrow s_{1,2}
=0,\pi $, we have in the first order in $s$, for $s\simeq 0,\pi $:
$$
{\hat {Q}}({\hat {d}}_i,s\simeq 0)\simeq 1+{{is} \over 2}{\hat {d}}_i 
$$
\begin{eqnarray}
{\hat {Q}}({\hat {d}}_i,s\simeq \pi)\simeq 1-{{i(\pi -s)} \over 2}{\hat {d}}_i.
\end{eqnarray}
By using eqs.(27,30) and $q^{{\hat {d}}_i} \simeq 1+is{\hat {d}}_i$, we have in the same
limiting cases: 
$${\hat {\partial }}_{x}\simeq {\partial }_{x}
\pm i\epsilon ({1 \over 2}{x\partial }_{x}^2+y{\partial }_{x}{\partial }_{y}+
z{\partial }_{x}{\partial }_{z}),
$$
$${\hat {\partial }}_{y}\simeq {\partial }_{y}
\pm i\epsilon ({1 \over 2}y{\partial }_{y}^2+z{\partial }_{z}{\partial }_{y}),
$$
\begin{eqnarray}
{\hat {\partial }}_{z}&\simeq &{\partial }_{z}
\pm i\epsilon {1 \over 2}z{\partial }_{z}^{2}. \hfill
\end{eqnarray}
where the first sign holds for $\epsilon =s \simeq 0$ and the second sign holds
for $\epsilon =\pi -s\simeq 0$. 

In the following we analyse a simple effect of non-commutativity
of the coordinate space on the non-relativistic dynamics of a quantum
particle. 
In order to get some physical information about such a
model,
in comparison with the normal plane, we  investigate some properties of
the operators of momentum and of the $z$ component of the angular momentum,
defined by ${\hat {L}}_{z}=
({\vec {r}}\times {\hat {{\vec {p}}}})_{z}$. 
The quantum non-commutative plane, with its non-commutative differential
structure is defined from eqs.(26-28) 
by the commutation relations:
\begin{eqnarray}
xy=qyx, \ \ \ \ {\hat {p}}_{x}=-iq^{2}\partial _{x}, \ \ \ \ {\hat {p}}_{y}=-iq\partial _{y},
\end{eqnarray}
and hence
$$
{\hat {p}}_{x}y=qy{\hat {p}}_{x}, \ \ \ \ {\hat {p}}_{y}x=qx{\hat {p}}_{y}, \ \ \ \ {\hat {p}}_{y}{\hat {p}}_{x}=q{\hat {p}}_{x}{\hat {p}}_{y}
$$
\begin{eqnarray}
{\hat {p}}_{x}x=-iq^{2}+q^{2}x{\hat {p}}_{x}+q(q-1)y{\hat {p}}_{y}, \ \ \ 
{\hat {p}}_{y}y=-iq+q^{2}y{\hat {p}}_{y}.
\end{eqnarray}
In the commutative plane, the generators:
${\hat {P}}_x =\partial _{x}$, ${\hat {P}}_{y}=\partial _{y}$ (translations or $-i$ times the 
momentum operators),  and ${\hat {R}}=y\partial _{x}-x\partial _{y}$ (rotation or 
the ${\hat {L}}_{z}$ component of the angular momentum operator) 
fulfil the commutator relations: 
\begin{eqnarray}
[{\hat {R}},{\hat {P}}_{x}]={\hat {P}}_{y}, \ \ \ \ [{\hat {R}},{\hat {P}}_{y}]=-{\hat {P}}_{x}, \ \ \ \ 
[{\hat {P}}_{x},{\hat {P}}_{y}]=0, 
\end{eqnarray}
and result in  a differential realisation of the Euclidean 
(Lie) algebra $E(2)$, of translations and rotations in the plane.
From the quantum mechanical point of view, (${\hat {P}}_{x,y},
{\hat {R}} \rightarrow 
{\hat {p}}_{x,y}, {\hat {L}}_{z}$), on the 
Hilbert space of states in the $(x,y)$ representation, we have an
uncertainty relation between ${\hat {p}}_{y}$ and ${\hat {L}}_{z}$ in the form
$4<{\hat {p}}_{y}^{2}><{\hat {L}}_{z}^{2}>\geq <{\hat {p}}_{x}^{2}>$.
In the non-commutative plane, eqs.(34)
do not close under the commutator relations and
we have not a closed q-algebra, like the case
of $E_{q}(2)$:
$$
[{\hat {p}}_{x},{\hat {p}}_{y}]=(1-q){\hat {p}}_{x}{\hat {p}}_{y},
$$
\begin{eqnarray}
[{\hat {p}}_{x},{\hat {L}}_{z}]=iq{\hat {p}}_{y}+(q-1){\hat {L}}_{z}
{\hat {p}}_{x}-q(q^{2}-1)y{\hat {p}}_{y}^{2},
\end{eqnarray}
$$
[{\hat {p}}_{y},{\hat {L}}_{z}]=-iq{\hat {p}}_{x}-(q^{3}-1){\hat {L}}_{z}{\hat {p}}_{y}+q(q^{2}-1)x{\hat {p}}_{y}^{2},
$$
where ${\hat {L}}_{z}=-iq(qy\partial _{x}-x\partial _{y})$.
Eqs.(35)
reduce to eqs.(34) when $q \rightarrow 1$ .
These commutator relations are quadratic and, in order to close this
q-deformed algebra, one needs to introduce the coordinate operators, too.
When $q\rightarrow 1$ eqs.(35)
become:
\begin{eqnarray}
[{\hat {p}}_{x},{\hat {p}}_{y}]=2{\hat {p}}_{x}{\hat {p}}_{y}, \ \ 
[{\hat {p}}_{x},{\hat {L}}_{z}]=-i{\hat {p}}_{y}-2{\hat {L}}_{z}{\hat {p}}_{x}, \ \ 
[{\hat {p}}_{y},{\hat {L}}_{z}]=i{\hat {p}}_{x}+2{\hat {L}}_{z}{\hat {p}}_{y},
\end{eqnarray}
and the q-algebra closes. It is a quadratic deformation of the algebra 
eqs.(34).
In the commutative case,
from the second commutator relation in eqs.(34),
the RHS is zero on a subspace of
the $(x,y)-$representation of the Hilbert space, given by wave functions
$\Psi (x,y)=\Psi (y)$. In this case both ${\hat {L}}_{z}$ and ${\hat {p}}_{y}$,
due to the above uncertainty relation, can be measured
with the same precision, i.e. there is no macroscopic motion in the $y$
direction. The wave function is a constant with respect to $x$ which forbids
it to belong to ${L}_{2}(R)$. We have full delocalisation in the $x$
direction and ${\hat {p}}_{x}\Psi =0$. 
In the non-commutative case we
look for such subspaces, wich annihilate the RHS in the
third commutator
relation in eqs.(35) and we get a nontrivial 
partial differential equation for $\Psi (x,y)$ 
$$
\biggl ( -q\partial _{x}+q(q^3 -1)y\partial _{x}\partial _{y}-
(q^2 -1)x(\partial _{x}^{2}+\partial _{y}^{2})-q^2 (q-1)
{\partial }_{y}^{2} \biggr ) \Psi
(x,y)=0,
$$
which, in the limit $q\rightarrow -1$ becomes
\begin{eqnarray}
(\partial _{x}+2y\partial _{x}\partial _{y}+2x{\partial }_{y}^{2})\Psi
(x,y)=0.
\end{eqnarray}
We search  solutions in the form $\Psi (x,y)=g(x)f(y)$. By introduction 
this form in eq.(37), performing the derivations by taking care of the
order of the operation in $x$ and $y$, 
we obtain one bounded $L_{2}(R)$ exact solution, 
in the form 
\begin{eqnarray}
\Psi (x,y)=e^{-\alpha x^2}\sqrt {y}{I}_{1/4}\biggl ( {{\alpha y^2 } \over {2}}
\biggr )
e^{-{{\alpha y^2} \over 2}},
\end{eqnarray}
where ${I}_{1/4}$ is the Bessel function of imaginary argument and
$\alpha $ is an arbitrary real parameter.
This solution represents a bounded function at $\infty $
but has one pole at $y=0$. Its asympthotic behaviour for $y \rightarrow
\infty $ is given by:
$\Psi (x,y)\simeq e^{-\alpha x^2}$ and describes a wavefunction which 
is constant with respect to $y$.
The wave function
is localised in the $x$ direction and
$<{\hat {L}}_{z}>|_{\Psi }=<{\hat {p}}_{x}>| _{\Psi }=<{\hat {p}}_{y}>| _{\Psi}=0$.
Consequently, the non-commutative plane provides 
a behaviour for the free particle similar with the existence of a
 potential valley in the $x$ direction. 

We note that in the general case, for an arbitrary $q\neq \pm 1$, 
eq.(37) becomes
a finite-difference {\bf {-}} partial differential equation, since, 
due to
the order dependent operations on $x,y$, the unknown function will appear in
these equations in the forms $\Psi (x,y)$ and $\Psi (x,qy)$, too.  

We are interested to find out
what implications result from the
use of the two different systems of coordinates $(x_{i},t)$ and
$({\hat x}_{i},t)$, the first being commutative
and the second non-commutative. 
We take a free 
non-relativistic particle described by
the Schr\"odinger equation, which
in the coordinates associated with the non-commutative space reads:
\begin{eqnarray}
i\hbar {{\partial \Psi } \over {\partial t}}
={1 \over {2m}}\biggl ( {\hat p}_x^2 +{\hat p}_y^2 +{\hat p}_z^2
\biggr ),
\end{eqnarray}
where ${\hat p}_{x_i}=-i\hbar {\partial }_{i}$ are the non-commutative
momentum operators defined in eqs.(27).
In the limits $s\rightarrow 0,\pi $
we are looking for solutions of eq.(39) in the perturbative form
$\Psi ({\vec r},t)=(\Phi _0+\epsilon \Phi 
({\vec r},t))e^{i({\vec k}{\vec r}-\omega t)}$
where $\Phi _0$ is a constant and $\epsilon =0$ for $s\simeq 0$
and $\epsilon =\pi -s$ for $s\simeq \pi $.
By taking into account eqs.(30,31) for $s\simeq 0$ and $\pi $,
we calculate, in first order in $\epsilon $ 
the approximate form of eq.(39), 
expressed back in the commutative coordinates. In this case
we introduce
a shift in the momentum operator of the form
\begin{eqnarray}
{\hat {\partial }}_{i}&=&{\partial }_{i}+A_{i}, \hfill
\end{eqnarray}
with
\begin{eqnarray}
{\vec {A}}=\mp s {\hbar} \biggl (
{{k_{x}^{2}} \over 2}x+k_{x}k_{y}y
+k_{x}k_{z}z,
{{k_{y}^{2}} \over 2}y+k_{y}k_{z}z,
{{k_{z}^{2}} \over 2}z \biggr ),
\end{eqnarray}
obtained from the action of the operator ${\hat {\partial }}_{i}$
from eq.(31) on the exponential function, i.e.
$
{\hat {\partial }}_{x}exp(i{\vec {k}}{\vec {r}})\simeq 
(ik_x \mp is ({1 \over 2}xk_{x}^{2}+y
k_x k_y +zk_x k_z))exp(i{\vec {k}}{\vec {r}})
$, etc.
The transformation of the coordinates and momenta, eqs.(3,24,27),
 results in the occurence of a
non-integrable phase in the wave function
\begin{eqnarray}
\Psi&=&{\Psi }^{'}exp \biggl ( i\int A_{i}dx_{i} \biggr ), \hfill
\end{eqnarray}
with $curl {\vec A}=\mp s{\hbar }(-k_{y}k_{z},k_{x}k_{z},k_{x}k_{y})\neq 0$.
If the motion is confined into a plane ($x_{i},x_{j}$) than the vector $curl
{\vec A}$ is always orthogonal to this plane.
A consequence of this approach is  the occurence of an anisotropy of the
wavefunction with respect to the coordinates,
resulting from the anysotropy of the non-commutative space, eq.(24). 
Consequently,
the non-integrable phase factor depends on the
direction of motion, i.e. on the wave function.
 
The similarity between the q-deformed momentum operators, and the kinetic
momentum operators:
\begin{eqnarray}
-i{\hbar }{\hat {\partial }}_{k}&\simeq &-i{\hbar }{\partial }_{k}-
{e \over c}A_{k},
\end{eqnarray}
is that one  of a  particle with charge $e$ moving in an
external magnetic field [7].
Using this analogy, we can say that the behaviour of a free particle
in a non-commutative system of coordinates, as described by
the Schr\"odinger equation, would be like that of a particle with charge $e$ moving
in a self-generated magnetic field given by
\begin{eqnarray}
B_{k}&\simeq &-{{{\hbar }c } \over e}sk_{j}k_{l},
\end{eqnarray}
with $(k,j,l)$ a permutation of $(1,2,3)$.
We also note that from eq.(41) results that this associated
magnetic field is not reducible to zero by performing 
a gauge transformation, i.e. 
$A_i$
can not be written in the form of a gradient.
We note here a possible connection between this approach and the anyon
statistics, [10]. The fractional statistics in the plane is 
obtained 
if the particles carry both charge and a magnetic flux. Due to the fact that
this flux is not involved in any restriction (neither integer nor
half-integer) such particles are anyons of fractional statistics.
In our case, i.e. the quantum plane, the essential feature is that fractional
statistics is implemented by means of the non-commutativity of the space.
In fact, rather than affixing by hand an internally prescripted vector
potential ${\vec {A}}$ to an ordinary particle, which is then transmuted
into an anyon, we claim that the apparition of ${\vec {A}}$
is the consequence of the non-commutativity of the space.
More, the corresponding Hamiltonian for a free particle in the non-commutative case,
expressed in terms of the commutative coordinates 
is similar with some typical examples of Hamiltonians of
the fractional quantum Hall effect [11].
However, in the present case the magnetic field is spreaded allover the plane,
different from the case of anyons and cyons, where the magnetic field is confined
only in a neighborhood of each particle.
Consequently, a different type of statistics occures here , since the
non-integrable phase-factor from eq.(42) depends on the contour of integration.
We can calculate the general exact solution for the Schr\"odinger equation
eq.(39) for the free particle in the non-commutative space with the
help of eqs.(27).
The solution of eq.(39) is given by a q-deformed exponential, instead
of the normal exponential, $\Psi ({\vec r} ,t)=
$ constant $e^{i({\vec k}{\vec r}-\omega t)}_q$. As $q\rightarrow 1$
this wavefunction reduces 
to the trivial plane waves of the commutative free space [8].
Another application of this result is the connection of this solution
with the solution for free electrons with quantum friction.
It is known that in this 
latter case for one-dimensional space, the wavefunctions have the form 
of a q-deformed exponential, too [8].

In conclusion, 
the idea of a smooth connection between the discrete and
continuous transformations, could lead to  new implications in the 
structure of space-time. First the space-time coordinates turn out to be
non-commutative
(of course in a sensible way only at very high energies) and invariant
at q-deformed groups of transformations.
The algebraic properties
with  physicalsemnifications remain only those invariant under such
nongeometric symmetries of the physical system.
The request of invariance under such symmetries
lead to modifications in the dynamics of 
quantum extended particles and of free
quantum particles. Loosely speaking, the broadening of the symmetries
from Lie to q-deformations
(in order to provide an appropriate algebraic frame for both the
 discrete and continuous transformations) could 
modify the statistics of identical particles in such 
spaces, and could restrict the freedom of motion of the particles due to
the unceratinty of observation of  the three coordinates, simultaneously.
\vfill
\eject

\vfill
\eject

{\bf {Figure captions}}
\vskip 1cm
Fig.1
\vskip 0.5cm
Real deformations of the
Coulomb potential $V_C (x)={{1} \over {x-1}}$, obtained with eq.(27), 
are plotted for different values of the deformation parameter
$q=e^s $, $s=0$, $-0.1$, $-0.25$, $-0.5$ and $-0.75$, on a logaritmic scale.
One can see that the pole of the q-deformed potential is translated in the
positive direction of the x-axis. 
\vskip 1cm
Fig.2
\vskip 0.5cm
Complex deformations of the same Coulomb potential as in Fig.1,
for different values of the deformation parameter
$q=e^{is} $: $s=0$ (the original $V(x)$ potential), 
$-0.15$, $-0.25$, $-0.5$, $-\pi /4$, $-\pi /2$,
$-\pi $ and $-10$ for the real part of $V(x,q)$, Fig.2a, and
$s=0$ (no imaginary part), $-0.15$, $-0.3$, $-9$ and
$-10$ for the imaginary part, Fig.2a.
The pole at $x=1$ is eliminated for $s\neq \pi Z$.
The values of the $-s$ parameter are shown next to each curve.
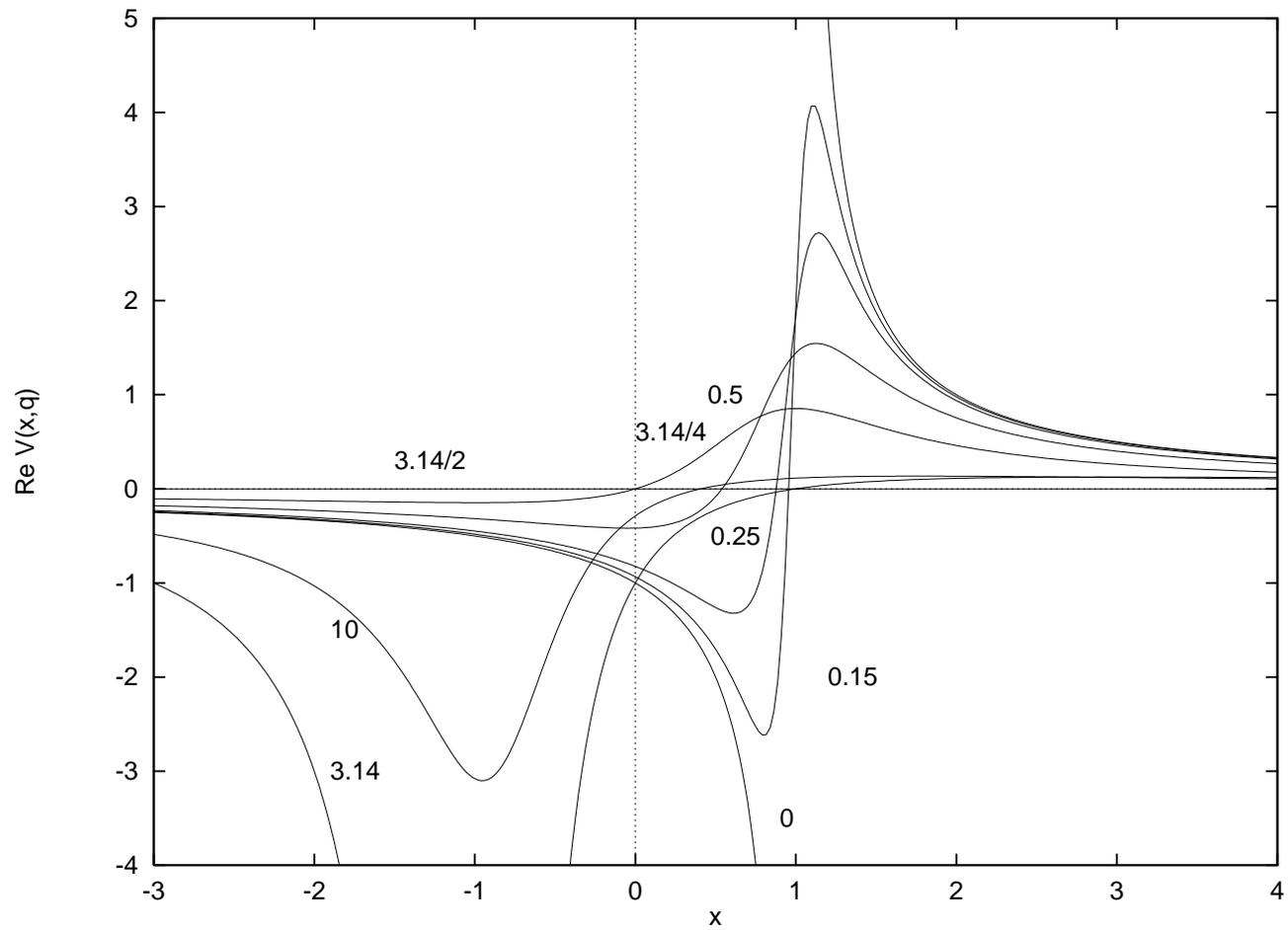
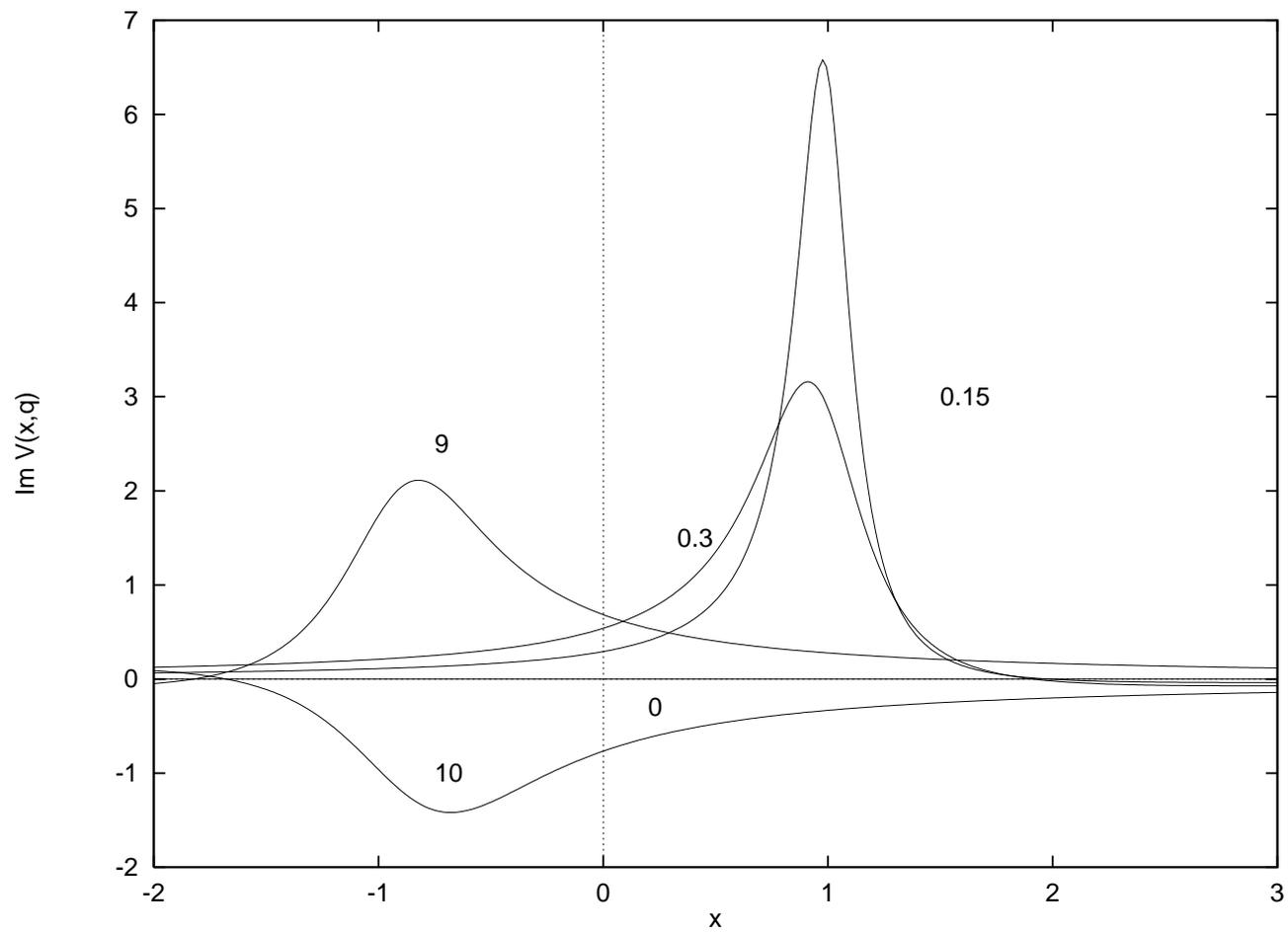
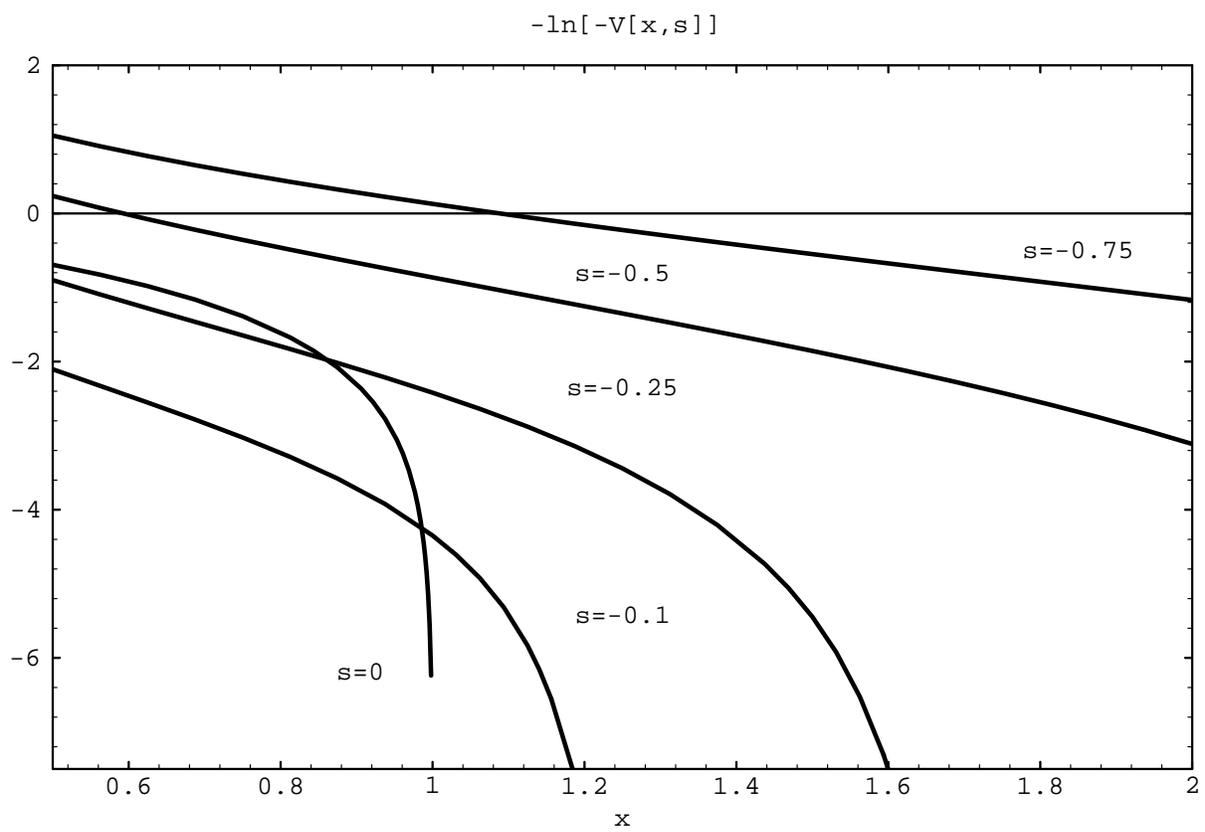
\end{document}